\documentclass[a4paper,12pt]{article}
\usepackage{enumerate}
\usepackage{latexsym}
\usepackage{amsfonts}
\usepackage[only,ninrm,elvrm,twlrm,sixrm,egtrm,tenrm]{rawfonts}
\usepackage{indentfirst}
\usepackage{amsmath}
\usepackage[noend]{algorithmic}
\usepackage{algorithm}
\usepackage{graphicx,psfrag}
\usepackage{graphics}
\usepackage{makeidx}

= 0.0in \topmargin = 0.0in \oddsidemargin=0.1in
\def\Dj{\hbox{D\kern-.73em\raise.30ex\hbox{-}
\raise-.30ex\hbox{}}}
\def\dj{\hbox{d\kern-.33em\raise.80ex\hbox{-}
\raise-.80ex\hbox{\kern-.40em}}}

\title{Partitioning 3-edge-colored complete\\ equi-bipartite graphs by
monochromatic trees\\ under a color degree
condition\footnote{Supported by NSFC, PCSIRT and the ``973"
program.}}
\author{\normalsize\itshape Xueliang Li and Fengxia Liu\\
\normalsize\itshape Center for Combinatorics and LPMC-TJKLC,\\
\normalsize\itshape Nankai University, Tianjin 300071, P. R. China\\
\normalsize lxl@nankai.edu.cn, xjulfx@163.com}
\date{}
\begin{document}
\maketitle
\begin{abstract}
The monochromatic tree partition number of an $r$-edge-colored graph
$G$, denoted by $t_r(G)$, is the minimum integer $k$ such that
whenever the edges of $G$ are colored with $r$ colors, the vertices
of $G$ can be covered by at most $k$ vertex-disjoint monochromatic
trees. In general, to  determine this number is very difficult. For
2-edge-colored complete multipartite graphs, Kaneko, Kano, and
Suzuki gave the exact value of $t_2(K(n_1,n_2,\cdots,n_k))$. In this
paper, we prove that if $n\geq 3$, and $K(n,n)$ is 3-edge-colored
such that every vertex has color degree 3, then $t_3(K(n,n))=3.$
\vspace{2mm}


\end{abstract}

\section{Introduction}

The {\it monochromatic tree partition number}, or simply {\it tree
partition number} of an $r$-edge-colored graph $G$, denoted by
$t_r(G)$, which was introduced by Erd\"{o}s, Gy\'{a}rf\'{a}s and
Pyber [1], is the minimum integer $k$ such that whenever the edges
of $G$ are colored with $r$ colors, the vertices of $G$ can be
covered by at most $k$ vertex-disjoint monochromatic trees.
Erd\"{o}s, Gy\'{a}rf\'{a}s and Pyber [1] conjectured that the tree
partition number of an $r$-edge-colored complete graph is $r-1$.
Moreover, they proved the conjecture for $r=3$ in the same paper.
For the case $r=2$, it is equivalent to the fact that for any graph
$G$, either $G$ or its complement is connected, an old remark of
Erd\"{o}s and Rado.

For infinite complete graph, Hajnal [2] proved that the tree
partition number for an $r$-edge-colored infinite complete graph is
at most $r$. For finite complete graph, Haxell and Kohayakawa [3]
proved that any $r$-edge-colored complete graph $K_n$ contains at
most $r$ monochromatic trees, all of different colors, whose vertex
sets partition the vertex set of $K_n$, provided $n\geq
3r^4r!(1-1/r)^{3(1-r)}\log r$. In general, to determine the exact
value of $t_r(G)$ is very difficult.

In this paper, we consider the tree partition number of complete
bipartite graphs. Notice that isolated vertices are also considered
as monochromatic trees. For any $m\geq n\geq 1$, let $K(A,B)=K(m,n)$
denote the complete bipartite graph with partite sets $A$ and $B$,
where $|A|=m$, $|B|=n$. Haxell and Kohayakawa [3] proved that the
tree partition number for an $r$-edge-colored complete bipartite
graph $K(n,n)$ is at most $2r$, provided $n$ is sufficiently large.
For 2-edge-colored complete multipartite graph
$K(n_1,n_2,\cdots,n_k)$, Kaneko, Kano, and Suzuki [5] proved the
following result: Let $n_1,n_2,\cdots n_k \ (k\geq 2)$ be integers
such that $1\leq n_1\leq n_2\leq \cdots\leq n_k$, and let
$n=n_1+n_2+\cdots +n_{k-1}$ and $m=n_k$. Then
$t_2(K(n_1,n_2,\cdots,n_k))=\lfloor \frac{m-2}{2^n}\rfloor +2$. In
particular, they proved that $t_2(K(m,n))=\lfloor
\frac{m-2}{2^n}\rfloor +2$, where $1\leq n\leq m$. Later in [4], Jin
et al gave a polynomial-time algorithm to partition a 2-edge-colored
complete multipartite graph into monochromatic trees.

In the present paper, we show that if $n\geq 3$ and $K(n,n)$ is
3-edge-colored such that every vertex has color degree 3, then
$t_3(K(n,n))=3$.

\section{Preliminaries}

In this section, we will give some notations and results on
2-edge-colored complete bipartite graphs. Although the result on the
partition number for 2-edge-colored complete bipartite graphs was
obtained by Kaneko, Kano and Suzuki in [5], and a polynomial-time
algorithm to get an optimal partition was obtained by Jin et al in
[4], in the following we will distinguish several cases, and for
each of which we will give the exact monochromatic trees to
partition the vertex set of a 2-edge-colored complete bipartite
graph. This gives not only the partition number for each case, but
more importantly, the clear structural description for the
partition, which will plays a key role for obtaining an optimal
partition in the 3-edge-colored case.

We first introduce two types of graphs. Let $G=K(A,B)$ be a
2-edge-colored complete bipartite graph, and all the edges are
colored with colors ``blue" or ``green". If the partite sets $A$ and
$B$ have partitions $A=A_1\cup A_2$ and $B=B_1\cup B_2$ with
$A_i\neq \emptyset$ and $B_i\neq \emptyset$ such that $K(A_1,B_1)$
and $K(A_2,B_2)$ are complete bipartite graphs with color blue,
$K(A_1,B_2)$ and $K(A_2,B_1)$ are complete bipartite graphs with
color green, then we call $K(A,B)$ an $M$-type graph. An $S$-type
graph is the graph satisfying $X(G)\neq \emptyset$ and $Y(G)\neq
\emptyset$, where $X(G)=\{u |$ all the edges incident with $u$ are
colored with blue $\}$, $Y(G)=\{u | \ \text{all the edges}$ incident
with $u$ are colored with green $\}$. Clearly, both $X(G)$ and
$Y(G)$ must be contained in a same partite set $A$ or $B$ of $G$. If
$G=K(A,B)$ is an $S$-type graph or an $M$-type graph, then we simply
denote it by $G\in S$ or $G\in M$.

In the following, if $G=K(A,B)$ is an $S$-type graph, and $X(G)\cup
Y(G)\in A$ (or, $X(G)\cup Y(G)\in B$), then we always denote
$A=A_1\cup A_2\cup A_3$ (or $B=B_1\cup B_2\cup B_3$), where
$A_1=\{u|$ all the edges incident with $u$ are colored blue$\}\
(B_1=\{u|$ all the edges incident with $u$ are colored blue$\})$,
$A_3=\{u|$ all the edges incident with $u$ are colored green$\}\
(B_3=\{u| $ all the edges incident with $u$ are colored green$\})$,
and $A_2=A-A_1-A_3\ (B_2=B-B_1-B_3)$. Clearly, $K(A_1\cup A_2,B)\
(K(B_1\cup B_2,A))$ has a blue spanning tree, and $K(A_3\cup A_2,B)\
(K(B_3\cup B_2,A))$ has a green spanning tree. That is, in the
following, we always use the subset with subscript ``1" to denote
the subset of vertices each of which is incident with only blue
edges, and the subset with subscript ``3" to denote the subset of
vertices each of which is incident with only green edges. For
example, if $K(R_2,C)$ is an $S$-type graph, and if we write
$R_2=R_{21}\cup R_{22}\cup R_{23}$, then $K(R_{21},C)$ is a maximum
blue complete bipartite subgraph of $K(R_2,C)$, and $K(R_{23},C)$ is
a maximum green complete bipartite subgraph of $K(R_2,C)$. If we
write $C=C_1\cup C_2\cup C_3$, then $K(R_2, C_1)$ is a maximum blue
complete bipartite subgraph of $K(R_2,C)$, and $K(R_2,C_3)$ is a
maximum green complete bipartite subgraph of $K(R_2,C)$.

\newtheorem{theorem}{Theorem}
\newtheorem{lemma}[theorem]{Lemma}
\begin{lemma} The 2-edge-colored complete bipartite graph $K(m,n)$ has
a monochromatic spanning tree if and only if $K(m,n)\notin S$ and
$K(m,n)\notin M$.
\end{lemma}
\noindent\textbf{Proof.} The necessity is obviously. Now we prove
the sufficiency.

If $K(m,n)$ has a vertex $x$ such that all the edges incident with
$x$ are colored blue (green), then every vertex is either in this
blue (green) star centered at $x$, or is adjacent to a vertex in the
blue (green) star by a blue (green) edge, since $K(m,n)\notin S$.
Thus, $K(m,n)$ has a blue (green) spanning tree.

We may assume therefore that for any vertex $x$ of $K(m,n)$, at
least one blue edge and one green edge are incident with it. Let $H$
be a subgraph of $K(m,n)$ induced by the green edges of $K(m,n)$,
and so $H$ is a spanning subgraph. If $H$ is connected, then $H$
contains a green spanning tree of $K(m,n)$, and the lemma follows.
Thus, we may assume that $H$ is not connected. Suppose $S$ is a
connected component of $H$, and $S\cap A=A_1$, $S\cap B=B_1.$ If
$A_1=A$, since every vertex in $B$ has at least one green edge
incident with it, we have $B_1=B$, which contradicts the assumption
that $H$ is disconnected. Thus, $A_1\neq A$ and $B_1\neq B$, and so
$K(A_1,B-B_1)$ and $K(A-A_1,B_1)$ are blue complete bipartite
graphs. Since $K(m,n)\notin M$, we have that at least one of
$K(A_1,B_1)$ and $K(A-A_1,B-B_1)$ is not green bipartite graph, and
so $K(A_1,B_1)$ and $K(A-A_1,B-B_1)$ have blue edges. Therefore,
$K(m,n)$ has a blue spanning tree. \hspace*{\fill}$\Box$

\vspace{3pt} Lemma 1 implies that if the 2-edge-colored complete
bipartite graph $K(m,n)$ does not have a monochromatic spanning
tree, then $K(m,n)\in S$ or $K(m,n)\in M$.

\begin{lemma} Let $K(A,B)$ be a 2-edge-colored complete bipartite graph.
If $K(A,B)\in M$, then the vertices of $K(A,B)$ can be covered by
two vertex-disjoint monochromatic trees with the same color.
\end{lemma}
\noindent\textbf{Proof.} Since $K(A,B)\in M$, we have partitions
$A=A_1\cup A_2$ and $B=B_1\cup B_2$ such that $K(A_1,B_1)$ and
$K(A_2,B_2)$ are blue complete bipartite graphs, $K(A_1,B_2)$ and
$K(A_2,B_1)$ are green complete bipartite graphs. That is, the
vertices of $K(A,B)$ can be covered by two vertex-disjoint blue
trees or two green trees.\hspace*{\fill}$\Box$

Let $K(A,B)$ be an $S$-type graph. If $A=A_1\cup A_2\cup
A_3(B=B_1\cup B_2\cup B_3)$, $|B|=1(|A|=1)$, and $|A_1|\geq
2(|B_1|\geq 2)$, $|A_3|\geq 2(|B_3|\geq 2)$, then we call $K(A,B)$
an $S_1^*$-type graph.

Let $K(A,B)$ be an $S$-type graph, and $A=A_1\cup A_2\cup A_3$. Then
for partition $B=B_i\cup \overline{B_i}$, we define

$b(B_i)=\{x\in A_2| \ K(x,B_i) \ \text{is a blue star, and} \
K(x,\overline{B_i})\  \text{is a green star}\}$,

$b(\overline{B_i})=\{x\in A_2| \ K(x,B_i) \ \text{is a green star,
and} \ K(x,\overline{B_i})\  \text{is a blue star}\}$.

If for any partition $B=B_i\cup \overline{B_i}$, $b(B_i)\neq
\emptyset$, $b(\overline{B_i})\neq \emptyset$, and $|A_1|\geq 2$,
$|A_3|\geq 2$, $|B|\geq 2$, then we call $K(A,B)$ an $S_1^{'}$-type
graph.

In the following, the $S_1^*$-type graphs and the $S_1^{'}$-type
graphs are denoted by $S_1$-type graph. The $S$-type graphs other
than the $S_1$-type graphs are denoted by $S_2$-type graph.

\begin{figure}[ht]
\begin{center}
\psfrag{3}{$A_1$} \psfrag{4}{$b(u_1)$} \psfrag{5}{$b(\{u_2,u_3\})$}
\psfrag{6}{$b(u_2)$} \psfrag{7}{$b(\{u_1,u_3\})$}
\psfrag{8}{$b(u_3)$} \psfrag{9}{$b(\{u_1,u_2\})$}
\psfrag{10}{$A_3$}\psfrag{11}{$u_1$}
\psfrag{12}{$u_2$}\psfrag{13}{$u_3$}
\includegraphics [width=12cm]{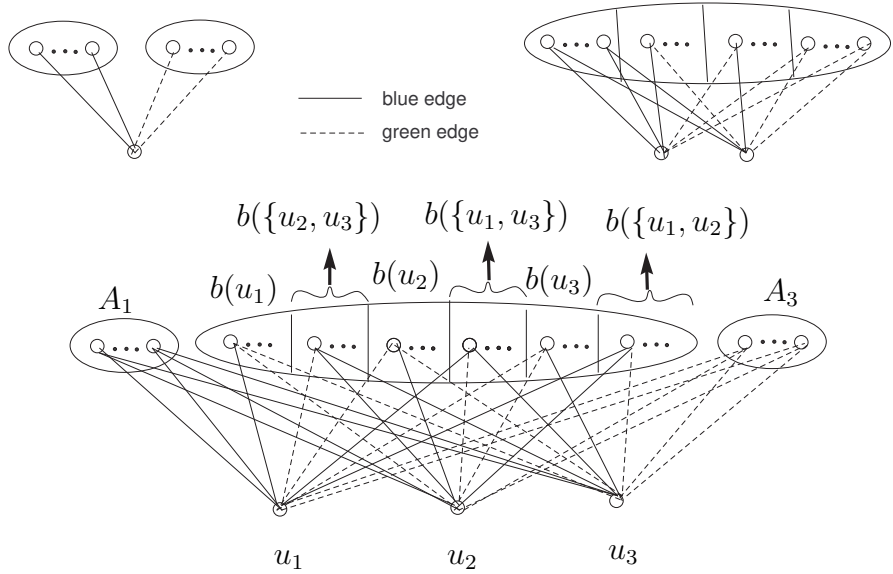}
\end{center}
\caption{$S_1$-type graphs.}
\end{figure}

Let $K(A,B)\in S$, and $A=A_1\cup A_2\cup A_3$. If $K(A,B)\in
S_1^*$, then $A_2=\emptyset$, and $|A|\geq 4=2^{|B|}+2$. If
$K(A,B)\in S_1^{'}$, then $A_2=\cup_{B=B_i\cup
\overline{B_i}}[b(B_i)\cup b(\overline{B_i})]$, here the union is
over all nonempty partitions of $B$, and for any $i$, $b(B_i)\neq
\emptyset$ and $b(\overline{B_i})\neq \emptyset$. Hence, $|A_2|\geq
2^{|B|}-2$, and so $|A|\geq 2^{|B|}+2$. Thus, if $K(A,B)\in S_1$,
then either $|A|\geq 2^{|B|}+2$ or $|B|\geq 2^{|A|}+2$ holds. If
$K(A,B)\in S_2$, either $min\{|A_1|,|A_3|\}=1$, or there exists a
partition $B=B_i\cup \overline{B_i}$ such that $b(B_i)=\emptyset$ or
$b(\overline{B_i})=\emptyset$.

\begin{lemma} Let $K(A,B)$ be a 2-edge-colored complete bipartite graph.
If $K(A,B)\in S_2$, then the vertices of $K(A,B)$ can be covered by
either an isolated vertex and a monochromatic tree or two
vertex-disjoint monochromatic trees with different colors.
Furthermore, if $A=A_1\cup A_2\cup A_3\ (B=B_1\cup B_2\cup B_3)$,
except the case $min\{|A_1|,|A_3|\}=1\ (min\{|B_1|,|B_3|\}=1)$, the
vertices of $K(A,B)$ always can be covered by two vertex-disjoint
monochromatic trees with different colors.
\end{lemma}
\noindent\textbf{Proof.} Without loss of generality, suppose
$A=A_1\cup A_2\cup A_3$.

\noindent \textbf{Case 1.} $min\{|A_1|,|A_3|\}=1$.

Since $K(A,B)\in S$, $K(A_1\cup A_2,B)$ and $K(A_3\cup A_2,B)$ have
a monochromatic spanning tree, respectively. Then the vertices of
$K(A,B)$ can be covered by an isolated vertex and a monochromatic
tree.

\noindent \textbf{Case 2.} There exists a partition $B=B_i\cup
\overline{B_i}$ such that $b(B_i)=\emptyset$ or
$b(\overline{B_i})=\emptyset$.

Without loss of generality, suppose $b(B_i)=\emptyset$. Let
$A_{21}=\{x\in A_2|\ K(x,\overline{B_i})$ have at least one blue
edge$\}$, and $A_{22}=A_{2}-A_{21}$. Since $b(B_i)=\emptyset$, every
vertex of $A_{22}$ has green edges to $B_i$. Then $K(A_1\cup
A_{21},\overline{B_i})$ has a blue spanning tree, and $K(A_3\cup
A_{22},B_i)$ has a green spanning tree. Thus, the vertices of
$K(A,B)$ can be covered by a blue tree and a green
tree.\hspace*{\fill}$\Box$

\begin{lemma} Let $K(A,B)$ be a 2-edge-colored complete bipartite graph.
Then $K(A,B)\in S_1$ if and only if the vertices of $K(A,B)$ can be
covered by at least three vertex-disjoint monochromatic trees.
\end{lemma}
\noindent\textbf{Proof.} We first consider the necessity. Without
loss of generality, suppose $A=A_1\cup A_2\cup A_3$. If $K(A,B)\in
S_1^*$, then $A_2=\emptyset$, $|B|=1$, and so the vertices of
$K(A,B)$ can be covered by at least $min\{|A_1|+1,|A_3|+1\}\geq 3$
vertex-disjoint monochromatic trees. For the case $K(A,B)\in
S_1^{'}$, if all the vertices of $B$ are in one monochromatic tree,
then the vertices of $K(A,B)$ can be covered by at least
$min\{|A_1|+1,|A_3|+1\}\geq 3$ vertex-disjoint monochromatic trees.
If all the vertices of $B$ are in two monochromatic trees, since for
any partition $B=B=B_i\cup \overline{B_i}$, we have $b(B_i)\neq
\emptyset$ and $b(\overline{B_i})\neq \emptyset$. So, the vertices
of $K(A,B)$ can be covered by at least $min_i\{min_{B=B_i\cup
\overline{B_i}}\{|b(B_i)|,|b(\overline{B_i})|\}+2\}\geq 3$
vertex-disjoint monochromatic trees. If all the vertices of $B$ are
in at least three monochromatic trees, then the vertices of $K(A,B)$
can be covered by at least three vertex-disjoint monochromatic
trees. In all cases, the vertices of $K(A,B)$ can be covered by at
least three vertex-disjoint monochromatic trees.

Now, we prove the sufficiency. If $K(A,B)\notin S_1$, then by the
above lemmas, the vertices of $K(A,B)$ can be covered by at most two
vertex-disjoint monochromatic trees, a contradiction.
\hspace*{\fill}$\Box$

From the above four lemmas, we have

\newtheorem{corollary}[theorem]{Corollary}
\begin{corollary} If $K(A,B)$ is a 2-edge-colored complete bipartite graph,
then it has one of the following four structures:

(1) $K(A,B)$ has a monochromatic spanning tree.

(2) $K(A,B)\in M$.

(3) $K(A,B)\in S_2$.

(4) $K(A,B)\in S_1$.
\end{corollary}

If $K(A,B)$ satisfies (2) or (3) of Corollary 5, then by Lemmas 2
and 3, the vertices of $K(A,B)$ can be covered by at most two
vertex-disjoint monochromatic trees. If $K(A,B)$ satisfies (4) of
Corollary 5, then from the proof of Lemma 4, the vertices of
$K(A,B)$ can be covered by
$min\{|A_1|+1,|A_3|+1,min_i|b(B_i)|+2,min_i|b(\overline{B_i})|+2\}$
vertex-disjoint monochromatic trees. Notice that
$min\{|A_1|+1,|A_3|+1,min_i|b(B_i)|+2,min_i|b(\overline{B_i})|+2\}\leq
\lfloor \frac{m-2}{2^n}\rfloor +2$, and the equality holds for some
graphs. So, the vertices of $K(A,B)$ can be covered by at most
$\lfloor \frac{m-2}{2^n}\rfloor +2$ vertex-disjoint monochromatic
trees, and there exists an edge coloring such that the vertices of
$K(A,B)$ are covered by exactly $\lfloor \frac{m-2}{2^n}\rfloor +2$
vertex-disjoint monochromatic trees. Thus, $t_2(K(m,n))=\lfloor
\frac{m-2}{2^n}\rfloor +2$.\vspace{5pt}

Let $K(A,B)$ be a 3-edge-colored complete bipartite graph, all the
edges of $K(A,B)$ are colored with colors ``red", ``blue" and
``green". Then we
define\\
\textbf{Case A.} All the vertices in $A$ are in some blue trees or
some green trees of $K(A,B)$.\\
\textbf{Case B.} All the vertices in $A$ are in some blue trees,
some green trees or a red tree of $K(A,B)$.\\
\textbf{Case C.} All the vertices in $A$ are in some blue trees,
some green trees or a set of isolated vertices of $K(A,B)$.

If we always consider the blue trees and the green trees first, and
the vertices in $B$ contained in these blue trees and green trees
are as small as possible, secondly, the red trees are as small as
possible, then the following lemma is obvious.

\begin{lemma} Let $K(A,B)$ be a 3-edge-colored complete bipartite
graph. If $|A|\leq|B|$, then there exists a tree partition
satisfying Case A or Case B. If $|A|>|B|$, then there exists a tree
partition satisfying Case A, Case B or Case C.
\end{lemma}

\section{Main result}

\begin{theorem} If $n\geq 3$, and $K(n,n)$ is 3-edge-colored such that
every vertex has color degree 3, then $t_3(k(n,n))=3.$
\end{theorem}
\noindent\textbf{Proof.} Assume that all the edges of
$K(n,n)(=K(A,B))$ are colored blue, green, or red. The vertices of
the graph in Figure 2 are covered by at least three vertex-disjoint
monochromatic trees. Then, $t_3(k(n,n))\geq 3$.

\begin{figure}[ht]
\begin{center}
\includegraphics[width=13cm]{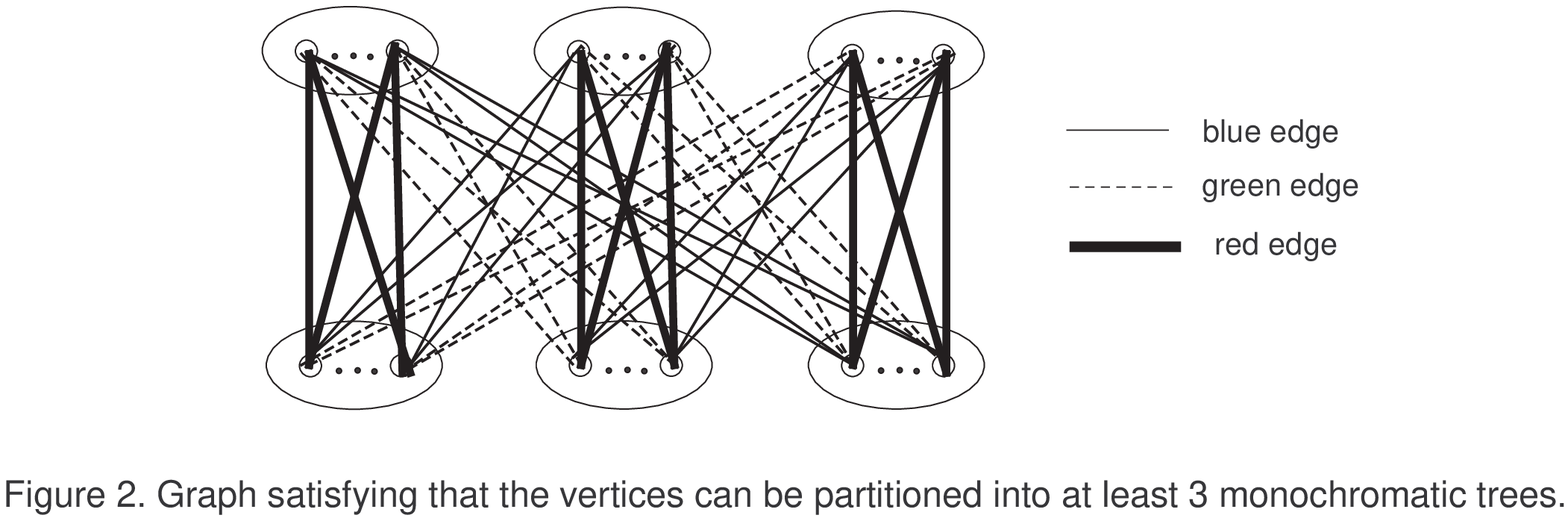}
\end{center}
\end{figure}

\vspace{-15pt}In the following, we prove $t_3(k(n,n))\leq 3$.
Suppose $R$ is the monochromatic connected component of $K(A,B)$
with the maximum number of vertices, without loss of generality,
suppose $R$ is red. Denote $R=R_1\cup R_2,\ R_1=R\cap A$, and
$R_2=R\cap B$.

If $R_1=A$, since the color degree of every vertex is 3, we have
$R_2=B$, then $K(A,B)$ has a red spanning tree.

We may assume therefore that $R_1\neq A$ and $R_2\neq B$. Denote
$C=A-R_1$ and $D=B-R_2$. Clearly, all the edges of $K(R_1,D)$ and
$K(R_2,C)$ are colored blue or green.

If the vertices of $K(C,D)$ can be covered by at most two
vertex-disjoint monochromatic trees, then the vertices of $K(A,B)$
can be covered by at most three vertex-disjoint monochromatic trees.
Thus, in the following, we assume that the vertices of $K(C,D)$ can
be covered by at least three vertex-disjoint monochromatic
trees.\vspace{5pt}

\noindent \textbf{Claim 1.} Every vertex in $K(C,D)$ has at least
one red edge incident with it, and there are at least one green edge
and one blue edge in $K(C,D)$.\vspace{5pt}

\noindent  \textbf{Proof.} Since every vertex has color degree 3,
and $K(R_1,D)$ and $K(R_2,C)$ are 2-edge-colored graphs with blue
and green, it is obvious that every vertex in $K(C,D)$ has at least
one red edge incident with it. Since the vertices of $K(C,D)$ can be
covered by at least three vertex-disjoint monochromatic trees, the
edges of $K(C,D)$ must be colored by at least two colors. Without
loss of generality, we assume $K(C,D)$ does not have green edges,
that is, $K(C,D)$ is a 2-edge-colored graph with blue and red. By
Lemma 4 we have $K(C,D)\in S_1$. Then, $K(C,D)$ has a vertex such
that all the edges incident with it are colored blue, which
contradicts the fact that every vertex in $K(C,D)$ has at least one
red edge incident with it. Thus, $K(C,D)$ has green edges.
\hspace*{\fill}$\Box$\vspace{5pt}

\noindent \textbf{Claim 2.} $|C|\geq 3$ and $|D|\geq 3$.\vspace{5pt}

\noindent \textbf{Proof.} Suppose $|C|\leq 2$. By Claim 1 every
vertex in $K(C,D)$ has at least one red edge incident with it, then
the vertices of $K(C,D)$ can be covered by two vertex-disjoint red
stars or a red spanning tree, which contradicts the assumption that
the vertices of $K(C,D)$ can be covered by at least three
vertex-disjoint monochromatic trees. \hspace*{\fill}$\Box$

Since $K(R_1,D)$ and $K(R_2,C)$ are 2-edge-colored graphs with blue
and green, by Corollary 5 we consider the following eight cases:\\

\noindent\textbf{Case 1.} Both $K(R_1,D)$ and $K(R_2,C)$ have
monochromatic spanning trees.

\noindent \textbf{Case 2.} One of $K(R_1,D)$ and $K(R_2,C)$ has a
monochromatic tree, the other is an $M$-type graph or an $S_2$-type
graph.

\noindent \textbf{Case 3.} One of $K(R_1,D)$ and $K(R_2,C)$ has a
monochromatic tree, the other is an $S_1$-type graph.

\noindent\textbf{Case 4.} $K(R_1,D)\in M$ and $K(R_2,C)\in M$.

\noindent\textbf{Case 5.} One of $K(R_1,D)$ and $K(R_2,C)$ is an
$M$-type graph, the other is an $S$-type graph.

\noindent\textbf{Case 6.} $K(R_1,D)\in S_2$ and $K(R_2,C)\in S_2$.

\noindent\textbf{Case 7.} $K(R_1,D)\in S_1$ and $K(R_2,C)\in S_1$.

\noindent\textbf{Case 8.} One of $K(R_1,D)$ and $K(R_2,C)$ is an
$S_1$-type graph, the other is an $S_2$-type graph.\\

In the following, we prove that for every above case, the vertices
of $K(A,B)$ can be covered by at most three vertex-disjoint
monochromatic trees.\vspace{5pt}

Clearly, in Case 1 the vertices of $K(A,B)$ can be covered by at
most two vertex-disjoint monochromatic trees. In Case 2, the
vertices of $K(A,B)$ can be covered by at most three vertex-disjoint
monochromatic trees.\\

\noindent\textbf{For Case 3,} without loss of generality, suppose
$K(R_1,D)$ has a green spanning tree, and $K(R_2,C)\in S_1$. Since
$K(R_2,C)\in S_1$, we have $|C|\geq 2^{|R_2|}+2$ or $|R_2|\geq
2^{|C|}+2$. Since $R$ is the maximum monochromatic component, and
$K(R_1,D)$ has a green spanning tree, we have $|D|\leq |R_2|$. If
$|C|\geq 2^{|R_2|}+2>2|R_2|$, then $|C|>2|R_2|=|R_2|+|R_2|\geq
|R_2|+|D|$, contradicting to $|R_1|+|C|=|R_2|+|D|=n$. If $|R_2|\geq
2^{|C|}+2$, then denote $R_2=R_{21}\cup R_{22}\cup R_{23}$. Since
every vertex has color degree 3, in $K(R_1,R_{21})$, every vertex in
$R_{21}$ has at least one green edge incident with it, and so
$K(R_1,R_{21}\cup D)$ has a green spanning tree. Obviously,
$K(C,R_{22}\cup R_{23})$ has a green spanning tree. Moreover, by
Claim 1, $K(C,D)$ has at least one green edge. Hence, $K(A,B)$ has a
green spanning tree, which contradicts our assumption that $R$ is
the maximum monochromatic component. Thus, this case does not
occur.\hspace*{\fill}$\Box$\\

\noindent\textbf{For Case 4,} we have $K(R_1,D)\in M$ and
$K(R_2,C)\in M$. By Lemma 2 the vertices of $K(R_1,D)$ and
$K(R_2,C)$ can be covered by two vertex-disjoint green trees,
respectively. By Claim 1 $K(C,D)$ has at least one green edge. Thus,
the vertices of $K(A,B)$ can be covered by at most three
vertex-disjoint green trees.\hspace*{\fill}$\Box$\\

\noindent\textbf{For Case 5,} without loss of generality, suppose
$K(R_1,D)\in M$, $K(R_2,C)\in S$. Since $K(R_2,C)\in S$, we have
$R_2=R_{21}\cup R_{22}\cup R_{23}$ or $C=C_1\cup C_2\cup C_3$. If
$R_2=R_{21}\cup R_{22}\cup R_{23}$, then $K(C,R_{22}\cup R_{23})$
has a green spanning tree. Since every vertex has color degree 3, in
$K(R_1,R_{21})$ every vertex in $R_{21}$ is incident with at least
one green edge. By Lemma 2 the vertices of $K(R_1,D)$ can be covered
by two vertex-disjoint green trees. Then, the vertices of
$K(R_1,R_{21}\cup D)$ can be covered by at most two vertex-disjoint
green trees. Moreover, $K(C,D)$ has at least one green edge. Thus,
the vertices of $K(A,B)$ can be covered by at most two
vertex-disjoint green trees. If $C=C_1\cup C_2\cup C_3$, by a
similar argument, the vertices of $K(R_1\cup C_1,D)$ can be covered
by at most two vertex-disjoint green trees, and $K(C_2\cup
C_3,R_{2})$ has a green spanning tree. Thus, the vertices of
$K(A,B)$ can be covered by at most three vertex-disjoint green
trees. \hspace*{\fill}$\Box$\\

\noindent\textbf{For Case 6,} we have $K(R_1,D)\in S_2$ and
$K(R_2,C)\in S_2$. Since $K(R_1,D)\in S_2$, we have $R_1=R_{11}\cup
R_{12}\cup R_{13}$ or $D=D_1\cup D_2\cup D_3$. Similarly, we have
$R_2=R_{21}\cup R_{22}\cup R_{23}$ or $C=C_1\cup C_2\cup
C_3$.\vspace{5pt}

\noindent\textbf{Subcase 6.1.} $R_1=R_{11}\cup R_{12}\cup R_{13}$
and $C=C_1\cup C_2\cup C_3$.\vspace{5pt}

Since every vertex has color degree 3, in $K(R_{11},R_{2})$ every
vertex in $R_{11}$ has at least one green edge incident with it. In
$K(C_1,D)$ every vertex in $C_1$ has at least one green edge
incident with it. Then $K(R_{11}\cup C_2\cup C_3,R_{2})$ and
$K(R_{12}\cup R_{13}\cup C_1,D)$ have a green spanning tree,
respectively. Thus, the vertices of $K(A,B)$ can be covered by at
most two vertex-disjoint green trees.\vspace{5pt}

\noindent\textbf{Subcase 6.2.} $R_2=R_{21}\cup R_{22}\cup R_{23}$
and $D=D_1\cup D_2\cup D_3$.\vspace{5pt}

The proof is similar to that of Subcase 6.1.\vspace{5pt}

\noindent\textbf{Subcase 6.3.} $R_1=R_{11}\cup R_{12}\cup R_{13}$
and $R_2=R_{21}\cup R_{22}\cup R_{23}$.\vspace{5pt}

Since $K(R_1,D)\in S_2$ and $K(R_2,C)\in S_2$, we can give the
following partition of $R_1,\ C,\ R_2$ and $D$, respectively:
$R_1=R_1^b\cup R_1^g$, $C=C_b\cup C_g$, $R_2=R_2^b\cup R_2^g$, and
$D=D_b\cup D_g$ such that $K(R_1^b,D_b)$ has a blue spanning tree,
$K(R_1^g,D_g)$ has a green spanning tree, $K(C^b,R_2^b)$ has a blue
spanning tree, and $K(C_g,R_2^g)$ has a green spanning tree.
Obviously, $R_{11}\subseteq R_1^b$, $R_{21}\subseteq R_2^b$,
$R_{13}\subseteq R_1^g$ and $R_{23}\subseteq R_2^g$. If
$K(R_{11},R_{21})$ has at least one blue edge, then $K(R_1^b,R_2^b)$
has at least one blue edge. Thus, the vertices of $K(A,B)$ can be
covered by at most one blue tree and two green trees. We may assume
therefore that all the edges of $K(R_{11},R_{21})$ are colored red
or green. Since $K(C,D)$ has at least one green edge,
$K(A-R_{11},B-R_{21})$ has a green spanning tree. If the vertices of
$K(R_{11},R_{21})$ can be covered by at most two vertex-disjoint
monochromatic trees, then the vertices of $K(A,B)$ can be covered by
at most three vertex-disjoint monochromatic trees. So, we assume
that the vertices of $K(R_{11},R_{21})$ can be covered by at least
three vertex-disjoint monochromatic trees. By Lemma 4 we have
$K(R_{11},R_{21})\in S_1$. Without loss of generality, we assume
that $R_{11}^r$ is the set with maximum number of vertices such that
$K(R_{11}^r,R_{21})$ is a red complete bipartite graph. Then
$K(R_{11}-R_{11}^r,R_{21})$ has a green spanning tree. Since every
vertex has color degree 3, in $K(R_{11}^r,R_{22}\cup R_{23})$ every
vertex in $R_{11}^r$ has at least one green edge incident with it,
and so $K(R_{11}^r\cup R_{12}\cup R_{13}\cup C,R_{22}\cup R_{23}\cup
D)$ has a green spanning tree. Thus, the vertices of $K(A,B)$ can be
covered by two vertex-disjoint green trees.\vspace{5pt}

\noindent\textbf{Subcase 6.4.} $C=C_1\cup C_2\cup C_3$ and
$D=D_1\cup D_2\cup D_3$.\vspace{5pt}

By the same arguments as in Case 6.3, we have partitions $C=C_b\cup
C_g$ and $D=D_b\cup D_g$. Clearly, $C_{1}\subseteq C_b$,
$C_{3}\subseteq C_g$, $D_{1}\subseteq D_b$ and $D_{3}\subseteq D_g$.
If $K(C_1,D_1)$ has at least one blue edge, or $K(C_3,D_3)$ has at
least one green edge, then the vertices of $K(A,B)$ can be covered
by at most three vertex-disjoint monochromatic trees. So, we assume
that $K(C_1,D_1)$ does not have blue edges, and $K(C_3,D_3)$ does
not have green edges. Then we have the following three subcases.
\vspace{5pt}

\noindent\textbf{Subcase 6.4.1.} $K(C_1,D_1)$ or $K(C_3,D_3)$ has a
monochromatic spanning tree.\vspace{5pt}

Since each of $K(C_2\cup C_3,R_2)$, $K(D_2\cup D_3,R_1)$, $K(C_1\cup
C_2,R_2)$ and $K(D_1\cup D_2,R_1)$ has a monochromatic spanning
tree, the vertices of $K(A,B)$ can be covered by at most three
vertex-disjoint monochromatic trees.\vspace{5pt}

\noindent\textbf{Subcase 6.4.2.} $K(C_1,D_1)\in S$ or $K(C_3,D_3)\in
S$.\vspace{5pt}

By a similar proof to the later part of Case 6.3, we can obtain that
the vertices of $K(A,B)$ can be covered by at most three
vertex-disjoint monochromatic trees.\vspace{5pt}

\noindent\textbf{Subcase 6.4.3.} $K(C_1,D_1)\in M$ and
$K(C_3,D_3)\in M$, see Figure 3.

\begin{figure}[ht]
\begin{center}
\psfrag{R_1}{$R_1$} \psfrag{R_2}{$R_2$} \psfrag{C_1}{$C_1$}
\psfrag{C_2}{$C_2$} \psfrag{C_3}{$C_3$} \psfrag{D_1}{$D_1$}
\psfrag{D_2}{$D_2$} \psfrag{D_3}{$D_3$}\psfrag{C_{11}}{$C_{11}$}
\psfrag{C_{12}}{$C_{12}$}\psfrag{C_{31}}{$C_{31}$}\psfrag{C_{32}}{$C_{32}$}
\psfrag{D_{11}}{$D_{11}$}\psfrag{D_{12}}{$D_{12}$}
\psfrag{D_{31}}{$D_{31}$}\psfrag{D_{32}}{$D_{32}$}
\includegraphics [width=12cm]{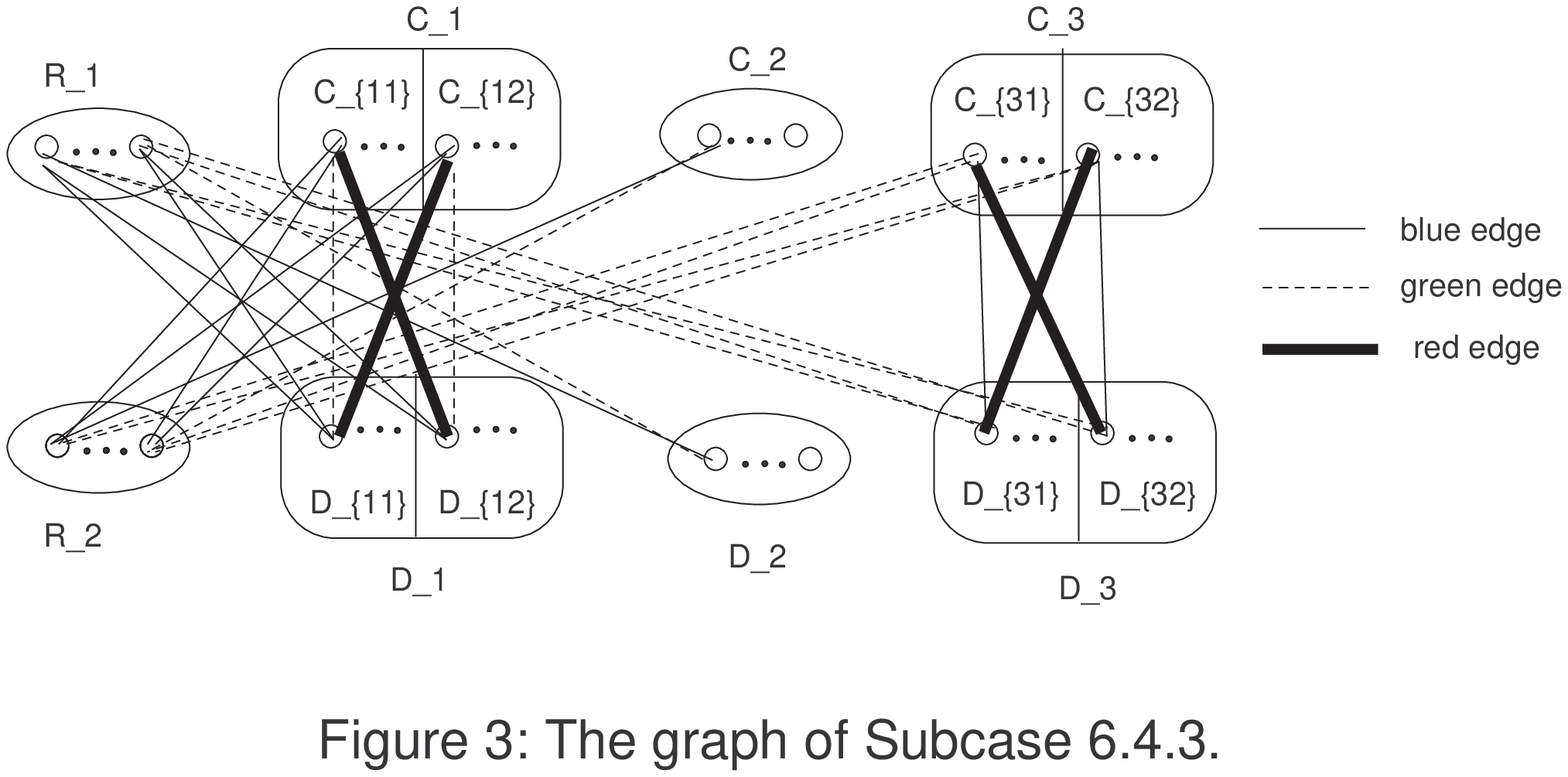}
\end{center}
\end{figure}

\vspace{-15pt}If all the edges of $K(C_1,D_3)$ are colored red, then
$K(C_{11}\cup C_3,D_3)$ has a red spanning tree. Since $K(C_{12}\cup
C_2,R_2)$ and $K(R_1,D_1\cup D_2)$ have blue spanning trees, the
vertices of $K(A,B)$ can be covered by three vertex-disjoint
monochromatic trees. Thus, we may assume that $K(C_1,D_3)$ has at
least one green edge or at least one blue edge. Without loss of
generality, assume $K(C_1,D_{31})$ has at least one blue edge. Then
$K(C_1\cup C_2\cup C_{31},R_2\cup D_{31})$, $K(C_{32},D_{32})$ and
$K(R_1,D_1\cup D_2)$ has blue spanning trees. Thus, the vertices of
$K(A,B)$ can be covered by three vertex-disjoint blue trees.\hspace*{\fill}$\Box$\\

\noindent\textbf{For Case 7,} we have $K(R_1,D)\in S_1$ and
$K(R_2,C)\in S_1$. Since $K(R_1,D)\in S_1$, we have $R_1=R_{11}\cup
R_{12}\cup R_{13}$ or $D=D_1\cup D_2\cup D_3$. Similarly, we have
$R_2=R_{21}\cup R_{22}\cup R_{23}$ or $C=C_1\cup C_2\cup
C_3$.\vspace{5pt}

\noindent\textbf{Subcase 7.1.} $R_1=R_{11}\cup R_{12}\cup R_{13}$
and $C=C_1\cup C_2\cup C_3$.\vspace{5pt}

Since $K(R_1,D)\in S_1$ and $K(R_2,C)\in S_1$, we have $|R_1|\geq
2^{|D|}+2>2|D|$ and $|C|\geq 2^{|R_2|}+2>2|R_2|$, and so
$|R_1|+|C|>2|R_2|+2|D|$, contradicting to $|R_1|+|C|=|R_2|+|D|=n$.
Thus, this case does not occur.\vspace{5pt}

\noindent\textbf{Subcase 7.2.} $R_2=R_{21}\cup R_{22}\cup R_{23}$
and $D=D_1\cup D_2\cup D_3$.\vspace{5pt}

The proof is similarly as Subcase 7.1.\vspace{5pt}

\noindent\textbf{Subcase 7.3.} $C=C_1\cup C_2\cup C_3$ and
$D=D_1\cup D_2\cup D_3$.\vspace{5pt}

Clearly, $K(C_1\cup C_2,R_2)$ and $K(C_2\cup C_3,R_2)$ have
monochromatic spanning trees, and $R$ is the maximum monochromatic
tree, then we have $|R_1|\geq \frac{|C|}{2}$. Similarly, $|R_2|\geq
\frac{|D|}{2}$. Moreover, by $K(R_2,C)\in S_1$, we have $|C|\geq
2^{|R_2|}+2>2|R_2|$. So, $|R_1|\geq \frac{|C|}{2}>|R_2|\geq
\frac{|D|}{2}$, that is, $|R_1|+|C|>|R_2|+|D|$, a contradiction.
Thus, this case does not occur.\vspace{5pt}

\noindent\textbf{Subcase 7.4.} $R_1=R_{11}\cup R_{12}\cup R_{13}$
and $R_2=R_{21}\cup R_{22}\cup R_{23}$.\vspace{5pt}

We define $R_{11}^{(2)}=\{u\in R_{11}|\ K(A,B)$ contains a green
$uv$-path for some $v\in R_1-R_{11}$ or $v\in R_2-R_{21}\}$,
$R_{11}^{(1)}=R_{11}-R_{11}^{(2)}$; $R_{13}^{(2)}=\{u\in R_{13}|\
K(A,B)$ contains a blue $uv$-path for some $v\in R_1-R_{13}$ or
$v\in R_2-R_{23}\}$, $R_{13}^{(1)}=R_{13}-R_{13}^{(2)}$;
$R_{21}^{(1)}, R_{21}^{(2)}, R_{23}^{(1)}$ and $R_{23}^{(2)}$ are
defined similarly.

Clearly, $K(R_{11}^{(1)},R_2-R_{21}^{(1)})$ and
$K(R_{21}^{(1)},R_1-R_{11}^{(1)})$ do not have green edges,
$K(R_{13}^{(1)},R_2-R_{23}^{(1)})$ and
$K(R_{23}^{(1)},R_1-R_{13}^{(1)})$ do not have blue edges. By Claim
2, $K(C,D)$ has at least one blue edge and one green edge, then
$K(A-R_{11}^{(1)},B-R_{21}^{(1)})$ has a green spanning tree,
$K(A-R_{13}^{(1)},B-R_{23}^{(1)})$ has a blue spanning tree. If the
vertices of $K(R_{11}^{(1)},R_{21}^{(1)})$ or
$K(R_{13}^{(1)},R_{23}^{(1)})$ can be covered by at most two
vertex-disjoint monochromatic trees, then the vertices of $K(A,B)$
can be covered by at most three vertex-disjoint monochromatic trees.
In the following, we consider the case that the vertices of both
$K(R_{11}^{(1)},R_{21}^{(1)})$ and $K(R_{13}^{(1)},R_{23}^{(1)})$
can be covered by at least three vertex-disjoint monochromatic
trees. We first give several remarks.\vspace{5pt}

\noindent\textbf{Remark 1.} $K(R_{11}^{(1)},R_{23}^{(1)})$ and
$K(R_{13}^{(1)},R_{21}^{(1)})$ are red complete bipartite graphs.

Since every vertex has color degree 3, we have \vspace{5pt}

\noindent\textbf{Remark 2.} Every vertex in
$K(R_{11}^{(1)},R_{21}^{(1)})$ has at least one green edge incident
with it, and every vertex in $K(R_{13}^{(1)},R_{23}^{(1)})$ has at
least one blue edge incident with it.

Since $R$ is the maximum monochromatic component, we
have\vspace{5pt}

\noindent\textbf{Remark 3.} $|R_{11}^{(1)}|+|R_{21}^{(1)}|\geq
|C|+|D|$ and $|R_{13}^{(1)}|+|R_{23}^{(1)}|\geq
|C|+|D|$.\vspace{5pt}

\noindent\textbf{Remark 4.} $\forall \ i=1,2,\ j=1,3, \
|R_{ij}^{(1)}|\geq 3$.\vspace{5pt}

\noindent\textbf{Proof.} Without loss of generality, suppose
$|R_{11}^{(1)}|\leq 2$. If $R_{11}^{(1)}=\emptyset$, since every
vertex has color degree 3, and $K(R_{21}^{(1)},A)$ has only blue
edges and red edges, we have $R_{21}^{(1)}=\emptyset$, which
contradicts to the assumption that the vertices of
$K(R_{11}^{(1)},R_{21}^{(1)})$ can be covered by at least three
vertex-disjoint monochromatic trees. If $1\leq |R_{11}^{(1)}|\leq
2$, then the vertices of $K(R_{11}^{(1)},R_{21}^{(1)})$ can be
covered by one green star or at most two vertex-disjoint green
trees, a contradiction. \hspace*{\fill}$\Box$\vspace{5pt}

\noindent \textbf{Remark 5.} $K(R_{11}^{(1)},R_{21}^{(1)})$ has at
least one red edge and one blue edge, $K(R_{13}^{(1)},R_{23}^{(1)})$
has at least one red edge and one green edge.\vspace{5pt}

\noindent\textbf{Proof.} If $K(R_{11}^{(1)},R_{21}^{(1)})$ does not
have red edges, by Remark 2, $K(R_{11}^{(1)},R_{21}^{(1)})$ does not
have any vertex such that all the edges incident with it are colored
blue, and so either $K(R_{11}^{(1)},R_{21}^{(1)})$ has a
monochromatic tree, or $K(R_{11}^{(1)},R_{21}^{(1)})\in M$, which
contradicts to the assumption that the vertices of
$K(R_{11}^{(1)},R_{21}^{(1)})$ can be covered by at least three
vertex-disjoint monochromatic trees. For other cases, we can prove
them similarly.\hspace*{\fill}$\Box$\vspace{5pt}

\noindent\textbf{Remark 6.} In $K(R_{11}^{(1)},R_{21}^{(1)}) \
(K(R_{13}^{(1)},R_{23}^{(1)}))$, every blue (green) edge has at
least one red edge and one blue (green) edge independent of it,
every red edge has at least one red edge and one blue (green) edge
independent of it.\vspace{5pt}

\noindent\textbf{Proof.} Let $e=uv$ be a blue edge of
$K(R_{11}^{(1)},R_{21}^{(1)})$. If $K(R_{11}^{(1)},R_{21}^{(1)})$
does not have red edges independent of $e$, then
$K(R_{11}^{(1)}-u,R_{21}^{(1)}-v)$ is a 2-edge-colored complete
bipartite graph with colors blue and green. If
$K(R_{11}^{(1)}-u,R_{21}^{(1)}-v)$ has a monochromatic spanning
tree, then the vertices of $K(R_{11}^{(1)},R_{21}^{(1)})$ can be
covered by at most two vertex-disjoint monochromatic trees, a
contradiction. If $K(R_{11}^{(1)}-u,R_{21}^{(1)}-v)\in M$, then the
vertices of $K(R_{11}^{(1)}-u,R_{21}^{(1)}-v)$ can be covered by two
vertex-disjoint green trees. Since $K(R_{11}^{(1)}-u,v)$ and
$K(R_{21}^{(1)}-v,u)$ both have green edges, the vertices of
$K(R_{11}^{(1)},R_{21}^{(1)})$ can be covered by at most two
vertex-disjoint green trees, a contradiction. If
$K(R_{11}^{(1)}-u,R_{21}^{(1)}-v)\in S$, noticing that
$K(R_{11}^{(1)}-u,v)$ and $K(R_{21}^{(1)}-v,u)$ have green edges,
then the vertices of $K(R_{11}^{(1)},R_{21}^{(1)})$ can be covered
by a green tree and a green star, a contradiction. Thus,
$K(R_{11}^{(1)},R_{21}^{(1)})$ has red edges independent of $e$. The
others can be proved similarly.\hspace*{\fill}$\Box$\vspace{5pt}

Since $|C|\geq 3$ and $|D|\geq 3$, we have $K(R_1,D)\in S_1^{'}$ and
$K(R_2,C)\in S_1^{'}$. Then $R_{12}=\cup_{D=D_i\cup
\overline{D_i}}[b(D_i)\cup b(\overline{D_i})]$ and
$R_{22}=\cup_{C=C_i\cup \overline{C_i}}[b(C_i)\cup
b(\overline{C_i})]$, here the union is over all nonempty partitions
of $D$ and $C$, respectively. For any nonempty partitions of $C$ and
$D$: $C=C_{i1}\cup \overline{C_{i2}}$, $D=D_{i1}\cup
\overline{D_{i2}}$, if $|b(C_{i1})|\geq |b(C_{i2})|$, then we denote
$C_{i1}=C_i$, $C_{i2}=\overline{C_i};$ if $|b(D_{i1})|\geq
|b(D_{i2})|$, then we denote $D_{i1}=D_i$, $D_{i2}=\overline{D_i}.$
So, in the following, if we write $C=C_i\cup \overline{C_i}$,
$D=D_i\cup \overline{D_i}$, then $|b(C_{i})|\geq
|b(\overline{C_i})|$ and $|b(D_{i})|\geq
|b(\overline{D_i})|$.\vspace{5pt}

\noindent\textbf{Subcase 7.4.1.} There exist partitions $C=C_k\cup
\overline{C_k}$ and $D=D_k\cup \overline{D_k}$ such that
$|b(C_{k})|\geq |b(\overline{D_k})|$ and $|b(D_{k})|\geq
|b(\overline{C_k})|$.\vspace{5pt}

In this case, $b(\overline{C_k})$ and $b(\overline{D_k})$ correspond
to the partite set $A$ in Lemma 4. Then by Lemma 6,
$K(b(D_{k}),b(\overline{C_k}))$ and $K(b(C_{k}),b(\overline{D_k}))$
have tree partitions satisfying Case A or Case B.\vspace{5pt}

\noindent\textbf{Subcase 7.4.1.1.} Both
$K(b(D_{k}),b(\overline{C_k}))$ and $K(b(C_{k}),b(\overline{D_k}))$
have tree partitions satisfying Case A.\vspace{5pt}

By Remark 5, $K(R_{11}^{(1)},R_{21}^{(1)})$ has at least one blue
edge, $K(R_{13}^{(1)},R_{23}^{(1)})$ has at least one green edge.
Then, $K(R_{11}\cup C_k,R_{21}\cup D_k)$ has a blue spanning tree,
and $K(R_{13}\cup \overline{C_k},R_{23}\cup \overline{D_k})$ has a
green spanning tree. By the definition of $b(D_{k})$ and $b(C_{k})$,
the vertices in $b(D_{k})$ and $b(C_{k})$ can be connected into the
blue tree of $K(R_{11}\cup C_k,R_{21}\cup D_k)$, and they also can
be connected into the green tree of $K(R_{13}\cup
\overline{C_k},R_{23}\cup \overline{D_k})$. Thus, the vertices of
$b(\overline{C_k})$ and $b(\overline{D_k})$ can be connected into
either the blue tree of $K(R_{11}\cup C_k,R_{21}\cup D_k)$ or the
green tree of $K(R_{13}\cup \overline{C_k},R_{23}\cup
\overline{D_k})$ by the vertices in $b(D_{k})$ and $b(C_{k})$.
Moreover, the vertices of $R_{12}-b(D_k)-b(\overline{D_k})$ have
either blue edges to $D_k$, or green edges to $\overline{D_k}$, the
vertices of $R_{22}-b(C_k)-b(\overline{C_k})$ have either blue edges
to $C_k$, or green edges to $\overline{C_k}$. Thus, the vertices of
$K(A,B)$ can be covered by a blue tree and a green tree.\vspace{5pt}

\noindent\textbf{Subcase 7.4.1.2.} One of
$K(b(D_{k}),b(\overline{C_k}))$ and $K(b(C_{k}),b(\overline{D_k}))$
has a tree partition satisfying Case A, the other has a tree
partition satisfying Case B.\vspace{5pt}

By a similar argument to that of Subcase 7.4.1.1, the vertices of
$K(A,B)$ can be covered by a blue tree, a green tree and a red
tree.\vspace{5pt}

\noindent\textbf{Subcase 7.4.1.3.} Both
$K(b(D_{k}),b(\overline{C_k}))$ and $K(b(C_{k}),b(\overline{D_k}))$
have tree partitions satisfying Case B.\vspace{5pt}

In $b(\overline{C_k})(b(\overline{D_k}))$, if the vertices in the
red tree satisfy that each of them has blue edges connecting to
$R_{11}(R_{21})$, or green edges connecting to $R_{13}(R_{23})$,
then all the vertices in $b(\overline{C_k})(b(\overline{D_k}))$ can
be connected into the blue tree of $K(R_{11}\cup C_k,R_{21}\cup
D_k)$ or the green tree of $K(R_{13}\cup \overline{C_k},R_{23}\cup
\overline{D_k})$ by the vertices in $b(D_{k})(b(C_k))$,
$R_{11}(R_{21})$ and $R_{13}(R_{23})$. Thus, the vertices of
$K(A,B)$ can be covered by a blue tree, a green tree and at most one
red tree. We may assume therefore that there exist vertices $x\in
b(\overline{C_k})$ and $y\in b(\overline{D_k})$ such that $x$ and
$y$ are the vertices in the red trees, and $K(x,R_{11}\cup R_{13})$,
$K(y,R_{21}\cup R_{23})$ are red stars. By Remark 5, we can suppose
that $uv$ is a red edge in $K(R_{11}^{(1)},R_{21}^{(1)})$, then the
red edges $ux$, $vy$ and $uv$ can connect the two red trees of
$K(b(D_{k}),b(\overline{C_k}))$ and $K(b(C_{k}),b(\overline{D_k}))$
into one red tree. By Remark 6, $K(R_{11}^{(1)}-u,R_{21}^{(1)}-v)$
has at least one blue edge. So, $K((R_{11}-u)\cup C_k,(R_{21}-v)\cup
D_k)$ still has a blue spanning tree. Thus, the vertices of $K(A,B)$
can be covered by a blue tree, a green tree and a red tree.
\vspace{5pt}

\noindent\textbf{Subcase 7.4.2.} For any partitions $C=C_i\cup
\overline{C_i}$ and $D=D_j\cup \overline{D_j}$, either
$|b(C_{i})|\geq |b(\overline{C_i})|>|b(D_{j})|\geq
|b(\overline{D_j})|$, or $|b(D_{j})|\geq
|b(\overline{D_j})|>|b(C_{i})|\geq |b(\overline{C_i})|$.\vspace{5pt}

Without loss of generality, suppose $C=C_k\cup \overline{C_k}$,
$D=D_k\cup \overline{D_k}$ such that $|b(D_{k})|\geq
|b(\overline{D_k})|>|b(C_{k})|\geq |b(\overline{C_k})|$. Define
$X_b(C_k)=\{x\in R_{22}|\ xu$ is a blue edge for some $u\in C_k$,
$xv$ is a green edge for some $v\in \overline{C_k}\}$,
$X_b(\overline{C_k})=\{x\in R_{22}|\ xu$ is a green edge for some
$u\in C_k$, $xv$ is a blue edge for some $v\in \overline{C_k}\}$.

Clearly, $b(C_{k})\subseteq X_b(C_k)$, $b(\overline{C_k})\subseteq
X_b(\overline{C_k})$ and $X_b(C_k)\cup X_b(\overline{C_k})=R_{22}$.
Then at least one of $|X_b(C_k)|\geq \frac{1}{2}|R_{22}|$ and
$|X_b(\overline{C_k})|\geq \frac{1}{2}|R_{22}|$ holds.\vspace{5pt}

\noindent\textbf{Subcase 7.4.2.1.} $|X_b(C_k)|\geq
\frac{1}{2}|R_{22}|$.\vspace{5pt}

In Subcase 7.4.1, we mainly use the property of $b(C_{k})$ that
every vertex in $b(C_{k})$ has blue edges to $C_k$ and has green
edges to $\overline{C_k}$. $X_b(C_k)$ also has the property. So, we
consider $K(b(D_{k}),b(\overline{C_k}))$ and
$K(b(\overline{D_k}),X_b(C_k))$ by the same argument as in Subcase
7.4.1. If $|X_b(C_k)|\geq |b(\overline{D_k})|$, then the vertices of
$K(A,B)$ can be covered by three vertex-disjoint monochromatic trees
just as Subcase 7.4.1. Hence, we consider the case $|X_b(C_k)|<
|b(\overline{D_k})|$. By Lemma 6, $K(b(D_{k}),b(\overline{C_k}))$
has a tree partition satisfying Case A or Case B, and
$K(b(\overline{D_k}),X_b(C_k))$ has a tree partition satisfying Case
A, Case B or Case C. If $K(b(\overline{D_k}),X_b(C_k))$ has a tree
partition satisfying Case A or Case B, then the proof is similar to
that of Subcase 7.4.1. If $K(b(\overline{D_k}),X_b(C_k))$ always has
a tree partition satisfying Case C, then denote the set of isolated
vertices in Case C as $I(\overline{D_k})$. If every vertex in
$I(\overline{D_k})$ has blue edges to $R_{21}$ or has green edges to
$R_{23}$, then similar to Subcase 7.4.1, the vertices of $K(A,B)$
can be covered by at most three vertex-disjoint monochromatic trees.

In the following, we assume that $I(\overline{D_k})$ has at least
one vertex such that all the edges incident with it in
$K(I(\overline{D_k}),R_{21}\cup R_{23})$ are colored red. Clearly,
$|b(\overline{D_k})-I(\overline{D_k})|\geq |X_b(C_k)|$, otherwise
$K(b(\overline{D_k}),X_b(C_k))$ has tree partition satisfying Case A
or Case B. So, $|b(\overline{D_k})|\geq
|I(\overline{D_k})|+|X_b(C_k)|$. Without loss of generality, suppose
$|R_{21}^{(1)}|\geq |R_{23}^{(1)}|$, then we consider
$K(I(\overline{D_k}),R_{21}^{(1)})$. Clearly, all the edges of
$K(I(\overline{D_k}),R_{21}^{(1)})$ are colored blue and
red.\vspace{5pt}

\noindent\textbf{Claim 3.} If $|I(\overline{D_k})|\leq
|R_{21}^{(1)}|$, then the vertices of $K(A,B)$ can be covered by at
most three vertex-disjoint monochromatic trees.\vspace{5pt}

\noindent\textbf{Proof.} In $K(I(\overline{D_k}),R_{21}^{(1)})$,
since $|I(\overline{D_k})|\leq |R_{21}^{(1)}|$, it is easy to see
that we have the fact that some vertices in $I(\overline{D_k})$ are
in blue trees and the others are in a red star. If
$K(b(D_{k}),b(\overline{C_k}))$ has a tree partition satisfying Case
A, or Case B such that in $b(\overline{C_k})$ all the vertices of
the red tree have green edges to $R_{13}$ or blue edges to $R_{11}$,
then the vertices of $K(A,B)$ can be covered by a blue tree, a green
tree and a red star. Otherwise, $K(b(D_{k}),b(\overline{C_k}))$
always has a tree partition satisfying Case B, and in
$b(\overline{C_k})$ there exists at least one vertex of the red tree
such that all the edges incident with it in
$K(b(\overline{C_k}),R_{13})$ are colored red. Then, similar to
Subcase 7.4.1.3, we can find a red edge $uv$ in
$K(R_{13}^{(1)},R_{23}^{(1)})$, and it can connect these two red
trees into one red tree, since $K(R_{21}^{(1)},R_{13}^{(1)})$ is a
red complete bipartite graph. Thus, the vertices of $K(A,B)$ can be
covered by a blue tree, a green tree and a red tree.
\hspace*{\fill}$\Box$\vspace{5pt}

If $|I(\overline{D_k})|> |R_{21}^{(1)}|$, then $|b(D_{k})|\geq
|b(\overline{D_k})|\geq
|I(\overline{D_k})|+|X_b(C_k)|>|R_{21}^{(1)}|+\frac{1}{2}|R_{22}|$,
and so
$|b(D_{k})|+|b(\overline{D_k})|>|R_{21}^{(1)}|+|R_{23}^{(1)}|+|R_{22}|$.
Thus, in Subcase 7.4.2.1, except
$|b(D_{k})|+|b(\overline{D_k})|>|R_{21}^{(1)}|+|R_{23}^{(1)}|+|R_{22}|$,
the vertices of $K(A,B)$ can be covered by at most three
vertex-disjoint monochromatic trees.\vspace{5pt}

\noindent\textbf{Subcase 7.4.2.2.} $|X_b(\overline{C_k})|\geq
\frac{1}{2}|R_{22}|$.\vspace{5pt}

In this case, we consider $K(b(D_{k}),b(C_{k}))$ and
$K(b(\overline{D_k}),X_b(\overline{C_k}))$. Since
$K(R_{11}^{(1)},R_{21}^{(1)})$ has at least one blue edge,
$K(R_{13}^{(1)},R_{23}^{(1)})$ has at least one green edge. We know
that $K(R_{11}\cup \overline{C_k},R_{21}\cup D_k)$ has a blue
spanning tree, and $K(R_{13}\cup C_k,R_{23}\cup \overline{D_k})$ has
a green spanning tree. We hope that the vertices of
$b(\overline{D_k})$ and $b(C_{k})$ can be connected to the blue tree
of $K(R_{11}\cup \overline{C_k},R_{21}\cup D_k)$ and the green tree
of $K(R_{13}\cup C_k,R_{23}\cup \overline{D_k})$, or they can
constitute a red tree. In this case, $b(C_k)$ and
$b(\overline{D_k})$ correspond to the partite set $A$ in Lemma 4.
Similar to Subcase 7.4.2.1, we can get the fact that except
$|b(D_{k})|+|b(\overline{D_k})|>|R_{21}^{(1)}|+|R_{23}^{(1)}|+|R_{22}|$,
the vertices of $K(A,B)$ can be covered by at most three
vertex-disjoint monochromatic trees.

By Subcase 7.4.2.1 and Subcase 7.4.2.2, we have that for partitions
$C=C_k\cup \overline{C_k}$ and $D=D_k\cup \overline{D_k}$ such that
$|b(D_{k})|\geq |b(\overline{D_k})|>|b(C_{k})|\geq
|b(\overline{C_k})|$, except
$|b(D_{k})|+|b(\overline{D_k})|>|R_{21}^{(1)}|+|R_{23}^{(1)}|+|R_{22}|$,
the vertices of $K(A,B)$ can be covered by at most three
vertex-disjoint monochromatic trees. If there exists a partition
$D=D_l\cup \overline{D_l}$ such that $|b(C_{k})|\geq
|b(\overline{C_k})|>|b(D_{l})|\geq |b(\overline{D_l})|$, then by a
similar argument to the above, we can obtain that except
$|b(C_{k})|+|b(\overline{C_k})|>|R_{11}^{(1)}|+|R_{13}^{(1)}|+|R_{12}|$,
the vertices of $K(A,B)$ can be covered by at most three
vertex-disjoint monochromatic trees. But
$|b(C_{k})|+|b(\overline{C_k})|>|R_{11}^{(1)}|+|R_{13}^{(1)}|+|R_{12}|$
contradicts to
$|b(C_{k})|+|b(\overline{C_k})|<|b(D_{k})|+|b(\overline{D_k})|<|R_{12}|$.
Thus, in the following we consider the case that for any partition
$D=D_i\cup \overline{D_i}$, we always have $|b(D_{i})|\geq
|b(\overline{D_i})|>|b(C_{k})|\geq |b(\overline{C_k})|$, and
$|b(D_{i})|+|b(\overline{D_i})|>|R_{21}^{(1)}|+|R_{23}^{(1)}|+|R_{22}|$,
otherwise, the vertices of $K(A,B)$ can be covered by at most three
vertex-disjoint monochromatic trees. Since $|D|\geq 3$, we have
$|R_{12}|=\sum_{D=D_i\cup \overline{D_i}}|b(D_i)\cup
b(\overline{D_i})|=\sum_{D=D_i\cup \overline{D_i}}[|b(D_i)|+|
b(\overline{D_i})|]$

\hspace{15pt}$>2(|R_{21}^{(1)}|+|R_{23}^{(1)}|+|R_{22}|)+|b(D_{k})|+|b(\overline{D_k})|$.\\
By Remark 3, $|R_{11}^{(1)}|+|R_{21}^{(1)}|\geq |C|+|D|$, that is,
$|R_{11}^{(1)}|-|D|\geq |C|-|R_{21}^{(1)}|$.\\
$|R_2|+|D|=|R_1|+|C|>|C|+|R_{11}|+|R_{13}|+2(|R_{21}^{(1)}|+|R_{23}^{(1)}|+|R_{22}|)$

\hspace{38pt}$+|b(D_{k})|+|b(\overline{D_k})|$.\\
$|R_{21}^{(2)}|+|R_{23}^{(2)}|>-|D|-|R_{22}|-|R_{21}^{(1)}|-|R_{23}^{(1)}|+
|C|+|R_{11}|+|R_{13}|$

\hspace{55pt}$+2(|R_{21}^{(1)}|+|R_{23}^{(1)}|+|R_{22}|)+|b(D_{k})|+|b(\overline{D_k})|$

\hspace{55pt}$\geq
2|C|-2|R_{21}^{(1)}|-|R_{22}|-|R_{23}^{(1)}|+|R_{11}^{(2)}|+|R_{13}|$

\hspace{55pt}$
+2(|R_{21}^{(1)}|+|R_{23}^{(1)}|+|R_{22}|)+|b(D_{k})|+|b(\overline{D_k})|$

\hspace{55pt}$>|b(D_{k})|+|b(\overline{D_k})|$.\\
Since $|b(D_{k})|\geq |b(\overline{D_k})|$, at least one of
$|R_{21}^{(2)}|>|b(\overline{D_k})|$ and
$|R_{23}^{(2)}|>|b(\overline{D_k})|$ holds. Without loss of
generality, we assume $|R_{21}^{(2)}|>|b(\overline{D_k})|$.

In the following, we consider $K(b(\overline{D_k}),R_{21}^{(2)})$
and $K(b(D_{k}),b(\overline{C_k}))$. $b(\overline{C_k})$ and
$b(\overline{D_k})$ correspond to the partite set $A$ in Lemma 4.
Clearly, $K(b(D_{k}),b(\overline{C_k}))$ has a tree partition
satisfying Case A or Case B. Let $X=\{\ v\in b(\overline{D_k})|\
K(v,R_{21}^{(2)})\text{ has blue edge}\}$, and $Y$ be the minimum
subset of $R_{21}^{(2)}$ satisfying that for any $v\in X$, there
exists a vertex $u\in Y$ such that $uv$ is a blue edge. Clearly,
$|Y|\leq |X|$. Denote $P=b(\overline{D_k})-X$ and
$Q=R_{21}^{(2)}-Y$. Then all the edges of $K(P,Q)$ are colored red
or green, and $|P|<|Q|$. We consider the following five small
cases.\vspace{5pt}

(1) $K(P,Q)$ has a green spanning tree.\vspace{5pt}

If $K(P,R_{23}^{(1)})$ has at least one green edge, then all the
vertices in $P$ and $Q$ can be connected to the green tree of
$K(R_{13}\cup \overline{C_k},R_{23}\cup \overline{D_k})$. So, the
vertices of $K(A,B)$ can be covered by at most three vertex-disjoint
monochromatic trees. We may assume therefore that
$K(P,R_{23}^{(1)})$ is a red complete bipartite graph. Let $uv$ be a
red edge in $K(R_{13}^{(1)},R_{23}^{(1)})$, then $K(P,v)$ is a red
star. Similarly, we can obtain that the vertices of $K(A,B)$ can be
covered by three vertex-disjoint monochromatic trees.\vspace{5pt}

(2) $K(P,Q)$ has a red spanning tree.\vspace{5pt}

If there exists a vertex $x\in P$ such that $K(x,R_{23}^{(1)})$ is a
red star, let $uv$ be a red edge in $K(R_{13}^{(1)},R_{23}^{(1)})$,
then $K(P,Q\cup v)$ has a red spanning tree. Similarly, the vertices
of $K(A,B)$ can be covered by three vertex-disjoint monochromatic
trees. Otherwise, every vertex in $P$ has green edge to
$R_{23}^{(1)}$, then it is easy to prove that the vertices of
$K(A,B)$ can be covered by at most three vertex-disjoint
monochromatic trees.\vspace{5pt}

(3) $K(P,Q)\in M$.\vspace{5pt}

Since $K(P,Q)\in M$, we can give the partitions $P=P_1\cup P_2$ and
$Q=Q_1\cup Q_2$ such that $K(P_1,Q_1)$ and $K(P_2,Q_2)$ are green
complete bipartite graphs, and $K(P_1,Q_2)$ and $K(P_2,Q_1)$ are red
complete bipartite graphs.

If both $K(P_1,R_{23}^{(1)})$ and $K(P_2,R_{23}^{(1)})$ have green
edges, then it is easy to prove that the vertices of $K(A,B)$ can be
covered by at most three vertex-disjoint monochromatic trees. If
both $K(P_1,R_{23}^{(1)})$ and $K(P_2,R_{23}^{(1)})$ do not have
green edges, then $K(P_1\cup P_2,R_{23}^{(1)})$ is a red complete
bipartite graph. Similarly, the vertices of $K(A,B)$ can be covered
by three vertex-disjoint monochromatic trees. Without loss of
generality, we may assume therefore that $K(P_1,R_{23}^{(1)})$ has a
green edge, say $wu$, and $K(P_2,R_{23}^{(1)})$ is a red complete
bipartite graph, then $K(P_2,R_{23}^{(1)}-u)$ is also a red complete
bipartite graph. Clearly, the vertices of $K(A,B)$ can be covered by
three vertex-disjoint monochromatic trees.\vspace{5pt}

(4) $K(P,Q)\in S_1$.\vspace{5pt}

Since $|P|<|Q|$, $K(P,Q)$ has a green tree containing all the
vertices in $P$. Then the proof is similar to the case that $K(P,Q)$
has a green spanning tree.\vspace{5pt}

(5) $K(P,Q)\in S_2$.\vspace{5pt}

Since $K(P,Q)\in S_2$, we can give partitions $P=P_r\cup P_g$ and
$Q=Q_r\cup Q_g$ such that $K(P_r,Q_r)$ has a red spanning tree,
$K(P_g,Q_g)$ has a green spanning tree. We consider four small
subcases.\vspace{5pt}

$\bullet$ At least one vertex in $P_g$ has green edge to
$R_{23}^{(1)}$, and every vertex in $P_r$ has green edge to
$R_{23}^{(1)}$.

$\bullet$ At least one vertex in $P_g$ that is incident with a green
edge to $R_{23}^{(1)}$, and at least one vertex in $P_r$ such that
all the edges incident with it are colored red.

$\bullet$ $K(P_g,R_{23}^{(1)})$ is a red complete bipartite graph,
and $K(P_r,R_{23}^{(1)})$ has at least one red edge.

$\bullet$ $K(P_g,R_{23}^{(1)})$ is a red complete bipartite graph,
and $K(P_r,R_{23}^{(1)})$ is a green complete bipartite
graph.\vspace{5pt}

For each of the above small subcases, we can easily obtain that the
vertices of $K(A,B)$ can be covered by two or three vertex-disjoint
monochromatic trees.\hspace*{\fill}$\Box$\\

\noindent\textbf{For Case 8,} without loss of generality, suppose
$K(R_1,D)\in S_1$, $K(R_2,C)\in S_2$. Since $K(R_1,D)\in S_1$, we
have $R_1=R_{11}\cup R_{12}\cup R_{13}$ or $D=D_1\cup D_2\cup D_3$.
Similarly, we have $R_2=R_{21}\cup R_{22}\cup R_{23}$ or $C=C_1\cup
C_2\cup C_3$. The case $R_1=R_{11}\cup R_{12}\cup R_{13}$ and
$C=C_1\cup C_2\cup C_3$, and the case $R_2=R_{21}\cup R_{22}\cup
R_{23}$ and $D=D_1\cup D_2\cup D_3$ are similar to Subcase 6.1 and
Subcase 6.2, respectively. The case $C=C_1\cup C_2\cup C_3$ and
$D=D_1\cup D_2\cup D_3$ is similar to Subcase 7.3. In the following,
we consider the case $R_1=R_{11}\cup R_{12}\cup R_{13}$ and
$R_2=R_{21}\cup R_{22}\cup R_{23}$. The proof is similar to that of
Subcase 7.4, by considering two subcases:\vspace{5pt}

\noindent\textbf{Subcase 8.1.} There exist partitions $C=C_k\cup
\overline{C_k}$ and $D=D_k\cup \overline{D_k}$ such that
$|b(C_{k})|\geq |b(\overline{D_k})|$ and $|b(D_{k})|\geq
|b(\overline{C_k})|$.\vspace{5pt}

This case can be proved similarly to Subcase 7.4.1.\vspace{5pt}

\noindent\textbf{Subcase 8.2.} For any partitions $C=C_i\cup
\overline{C_i}$ and $D=D_j\cup \overline{D_j}$, either
$|b(C_{i})|\geq |b(\overline{C_i})|>|b(D_{j})|\geq
|b(\overline{D_j})|$, or $|b(D_{j})|\geq
|b(\overline{D_j})|>|b(C_{i})|\geq |b(\overline{C_i})|$.\vspace{5pt}

Since $K(R_2,C)\in S_2$, we have $R_{22}=\cup_{C=C_i\cup
\overline{C_i}}[b(C_i)\cup b(\overline{C_i})]$ such that for some
$\overline{C_i}$, $b(\overline{C_i})=\emptyset$. Without loss of
generality, suppose $b(\overline{C_l})=\emptyset$, then we consider
partitions $C=C_l\cup \overline{C_l}$ and $D=D_j\cup \overline{D_j}$
for some $j$, hence $|b(D_{j})|\geq
|b(\overline{D_j})|>|b(C_{l})|\geq |b(\overline{C_l})|$. Thus, we
can prove it similarly to Subcase 7.4.2.\\

Up to now, we have exhausted all cases, and proved that for any
3-edge-colored complete bipartite graph $K(n,n)$ satisfying the
condition of Theorem 7, the vertices of it can be covered by at most
three vertex-disjoint monochromatic trees. Thus $t_3(k(n,n))\leq
3$.\hspace*{\fill}$\Box$

\section{Conclusion}

As one can see, we only considered 3-edge-colored complete bipartite
graphs with ``equal bipartition", and the ``color degree" of every
vertex is 3. These restrictions are really very helpful to
concluding our proofs. Even though, the proof looks very long and
complicated. More general questions are: can we drop the equal
bipartition restriction to get the partition number ? can we drop
the color degree restriction to get the partition number ? or can we
drop both restrictions to get the partition number ? We tried for a
year but failed to complete it. Things become out of control without
any of the restrictions.


\begin{thebibliography}{111}

\bibitem{1} P. Erd\"{o}s, A. Gy\'{a}rf\'{a}s and L. Pyber, Vertex coverings by
monochromatic cycles and trees, J. Combin. Theory, Ser. B 51 (1991),
90-95.

\bibitem{2} A. Hajnal, P. Komj\'{a}th, L. Soukup and I. Szalkai,
Decompositions of edge colored infinite complete graphs, Colloq.
Math. Soc. J\'{a}nos Bolyai 52 (1987), 277-280.

\bibitem{3} P.E. Haxell and Y. Kohayakawa, Partitioning by monochromatic trees, J. Combin.
Theory. Ser. B 68 (1996), 218-222.

\bibitem{4} Z.M. Jin, M. Kano, X. Li and B. Wei, Partitioning
2-edge-colored complete multipartite graphs into monochromatic
cycles, paths and trees, J. Comb. Optim 11(2006), 445-454.

\bibitem{5} A. Kaneko, M. Kano and K. Suzuki, Partitioning complete multipartite graphs by
monochromatic trees, J. Graph Theory 48 (2005), 133-141.

\end{thebibliography}
\end{document}